\documentclass[11pt,a4paper]{article}

\usepackage{bezier,amsfonts,amssymb,graphicx,amsthm}
\usepackage[english]{babel}
\def\lesta{ \hfill $\Box$ \bigskip}

\newcommand{\auttfg}{$\mbox{Aut}^{\mbox{{\tiny{\textbf{TF}}}}}G$}

\newcommand{\BG}{\mbox{\textbf{B}}(G)}

\newcommand{\aut}{\mbox{Aut}}

\newcommand{\id}{{\rm id}}

\newcommand{\ZZ}{\mathbb Z}

\newtheorem{Thm}{Theorem}[section]

\newtheorem{Prop}[Thm]{Proposition}

\newtheorem{Lem}[Thm]{Lemma}

\begin{document}
 \begin{titlepage}
\title{Unstable Graphs: \\ A Fresh Outlook via TF-Automorphisms}
\author{ J. Lauri, R. Mizzi \\ Department of Mathematics \\ University of Malta
\\ Malta \\ josef.lauri@um.edu.mt \\ russell.mizzi@um.edu.mt \and R. Scapellato \footnote{corresponding author} \\ Dipartimento 
di Matematica \\ Politecnico di Milano \\ Milano \\ Italy \\ raffaele.scapellato@polimi.it}
\pdfinfo{ /Author (Josef Lauri, Russell Mizzi, Raffalele Scapellato) /Title (Unstable Graphs: A Fresh Outlook) /Keywords (graphs;canonical double covers;two-fold isomorphisms) }\maketitle

\begin{small}
\begin{abstract}
In this paper, we first establish the very close link between stability of graphs, a concept first introduced in \cite{Scapsalvi1} and studied most notably by Surowski \cite{Surowski1}, \cite{Surowski2} and Wilson \cite{Wilson01} and two-fold automorphisms. The concept of two-fold isomorphisms, as far as we know, first appeared in literature in the form of  isotopies of digraphs \cite{zelinka4}, \cite{zelinka1}, \cite{zelinka2}, \cite{zelinka3} and later studied formally in \cite{lms1}, \cite{lms2} with a greater emphasis on undirected graphs.  We then turn our attention to the stability of graphs which have every edge on a triangle, but with the fresh outlook provided by TF-automorphisms. Amongst such graphs are strongly regular graphs with certain parameters. The advantages of this fresh outlook are highlighted when we ultimately present a method of constructing and generating unstable graphs with large diameter having every edge lying on a triangle. This was a rather surprising outcome. \end{abstract}
\end{small}
\bigskip

\bigskip 

\bigskip

{Keywords: \\
 \indent graph stability, canonical double covers, two-fold isomorphisms}

\end{titlepage}

 
 \section{General Introduction and Notation}
 
 Let $G$ and $H$ be simple graphs, that is, undirected and without loops or multiple edges. Consider the edge $\{u,v\}$ to be the set of arcs $\{(u,v),\ (v,u)\}$. A \textit{two-fold isomorphism} or TF-\textit{isomorphism} from $G$ to $H$ is a pair of bijections $\alpha$, $\beta$: V$(G)$ $\rightarrow$ V$(H)$ such that $(u,v)$ is an arc of $G$ if and only if it is an arc of $H$. When such a pair of bijections exist, we say that $G$ and $H$ are TF-\textit{isomorphic} and the TF-isomorphism is denoted by $(\alpha,\beta)$. The inverse of $(\alpha,\beta)$, that is, $(\alpha^{-1},\beta^{-1})$ is a TF-isomorphism from $H$ to $G$. Furthermore, if $(\alpha_{1},\beta_{1})$ and $(\alpha_{2},\beta_{2})$ are both TF-isomorphisms from $G$ to $H$ then so is $(\alpha_{1}\alpha_{2},\beta_{1}\beta_{2})$. When $\alpha=\beta$, the TF-isomorphism can be identified with the isomorphism $\alpha$.

\begin{figure}[h]
 \centering
 \includegraphics[width=8cm,height=4.6cm]{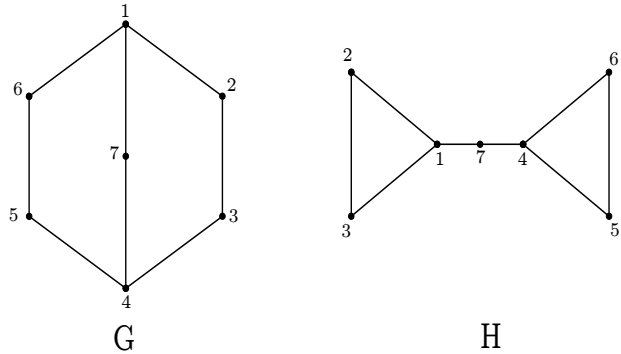}
 \caption{$G$ and $H$ are two non-isomorphic TF-isomorphic graphs. }\label{fig:tfisomorph01}
\end{figure}

The two graphs $G$, $H$ in Figure \ref{fig:tfisomorph01}, which have the same vertex set V$(G) =$ V$(H)$, are non-isomorphic and yet they are TF-isomorphic. In fact  $(\alpha,\beta)$ where $\alpha = (2\ 5)(1\ 4)(3\ 6)(7)$ and $\beta = (1\ 4)(3\ 6)$ is a TF-isomorphism from $G$ to $H$.\\

Some graph properties are preserved by a TF-isomorphism. Such is the case with the degree sequence, as illustrated by Figure \ref{fig:tfisomorph01}.  We also know that two graphs are TF-isomorphic if and only if that they have isomorphic canonical double covers \cite{lms1}.  Alternating paths or Z-trails, which we shall define in full below, are invariant under TF-isomorphism.  For instance, the alternating path $5 \longrightarrow 6 \longleftarrow 1 \longrightarrow 2$  in $G$  is mapped by $(\alpha,\beta)$ to the similarly alternating path $2 \longrightarrow 3 \longleftarrow 1 \longrightarrow 2$ which we shall later be calling "semi-closed".

\section{Notation}

 A \textit{mixed graph} is a pair $G=(\mbox{V}(G),\mbox{A}(G))$ where V$(G)$ is a set and A$(G)$ is a set of ordered pairs of elements of V$(G)$. The elements of V$(G)$ are called \textit{vertices} and the elements of A$(G)$ are called \textit{arcs}. When referring to an arc $(u,v)$, we say that $u$ is \textit{adjacent to} $v$ and $v$ is \textit{adjacent from} $u$. The ordered pair $(u,v)$ will often be denoted by $uv$. Sometimes we use $u \longrightarrow_{G}v$ or simply $u \longrightarrow v$ to represent an arc $(u,v) \in$ A$(G)$. The vertex $u$ is the \textit{start-vertex} and $v$ is the \textit{end-vertex} of a given arc $(u,v)$. An arc of the form $(u,u)$, or equivalently $uu$ is called a \textit{loop}.  A mixed graph cannot contain multiple arcs, that is, it cannot contain the arc $(u,v)$ more than once. A mixed graph $G$ is called \textit{bipartite} if there is a partition of V$(G)$ into two sets $X$ and $Y$, which we call \textit{colour classes}, such that for each arc $uv$ of G the set $\{u,v\}$ 
intersects both $X$ and $Y$. A set $S$ of arcs is \textit{self-paired} if, whenever $(u,v) \in$ $S$, $(v,u)$ is also in $S$. If $S$ $\ =\{(u,v), (v,u)\}$, then we consider $S$ to be the unordered pair $\{u,v\}$; this unordered pair is called an \textit{edge}.\\

It is useful to consider two special cases of mixed graphs. A \textit{graph} is a mixed graph without loops whose arc-set is self-paired. The edge set of a graph is denoted by E$(G)$. A \textit{digraph} is a mixed graph with no loops in which no set of arcs is self-paired. The \textit{inverse} $G'$ of a mixed graph $G$ is obtained from $G$ by reversing all its arcs, that is V$(G') =$V$(G)$ and $(v,u)$ is an arc of $G'$ if and only if $(u,v)$ is an arc of $G$. A digraph $G$ may therefore be characterised as a mixed graph for which A$(G)$ and A$(G')$ are disjoint. Given a mixed graph $G$ and a vertex $v \in$ V$(G)$, we define the \textit{in-neighbourhood} $N_{in}(v)$ by $N_{in}(v) = \{x \in \mbox{V}(G)| (x,v) \in \mbox{A}(G)\}$. Similarly we define  the  \textit{
out-neighbourhood} $N_{out}(v)$ by $N_{out}(v) = \{x \in \mbox{V}(G)| (v,x) \in \mbox{A}(G)\}$.  The \textit{in-degree} $ \rho_{in}(v)$ of a vertex $v$ is defined by $ \rho_{in}(v) = |N_{in}(v)|$ and the \textit{out-degree} $ \rho_{out}(v)$ of a vertex $v$ is defined by $ \rho_{out}(v) = |N_{out}(v)|$. When $G$ is a graph, these notions reduce to the usual neighbourhood $N(v)=N_{in}(v)=N_{out}(v)$ and degree $\rho(v)=\rho_{in}(v)=\rho_{out}(v)$. \\

Let $G$ be a graph and let $v \in$ V$(G)$. Let $N(v)$ be the neighbourhood of $v$. We say that $G$ is \emph{vertex-determining} if $N(x)\neq N(y)$ for any two distinct vertices $x$ and $y$ of $G$ \cite{Scapsalvi1}.\\

A set $P$ of arcs is called an \textit{alternating trail} or Z-\textit{trail} if its elements can be ordered in a sequence $a_1, a_2, \cdots,a_k$ such that if $a_{i}=(u,v)$ or $a_{i+1}=(u,w)$,  $w\neq v$ or  $a_{i+1}=(w,v)$, $w\neq u$.  Note that in a Z-trail vertices may be repeated but arcs may not.  If $a_{1} =(p,q)$ and $a_{k}=(r,s)$ we then say that the Z-trail joins $p$ and $s$.  A Z-trail can be \textit{open} (when the first and last vertex are different), \textit{semi-closed} (when the first arc is $uv$ and the last arc is $wu$) or \textit{closed} (when the first arc is $uv$ and the last arc is  $uw$ or the first arc is $vu$ and the last arc is $wu$).  Note that the ordering of arcs described in the definition is not unique, but it is basically unique for open and semi-closed Z-trails (one just chooses the beginning or the end), while for closed Z-trails all vertices can be taken as the first vertex. Whether a Z-trail is open, semi-closed or closed is determined by the first and last vertices. Open trails can have any number of arcs, semi-closed trails have an odd number of arcs and closed trails have an even number of arcs.\\


Any other graph theoretical terms which we use are standard and can be found in any graph theory textbook such as \cite{bondy}. For information on automorphism groups, the reader is referred to \cite{lauri2}.\\

Let $G$ and $H$ be two mixed graphs and suppose that $\alpha$, $\beta$ are bijections from V$(G)$ to V$(H)$. The pair $(\alpha,\beta)$ is said to be a \textit{two-fold isomorphism} (or TF\textit{-isomorphism}) if the following holds: $(u,v)$ is an arc of $G$ if and only if $(\alpha(u),\beta(v))$ is an arc of $H$. We then say that $G$ and $H$ are TF-\textit{isomorphic} and write $G\cong ^{\mbox{{\tiny{\textbf{TF}}}}} H$. Note that when $\alpha=\beta$ the pair $(\alpha,\beta)$ is a TF-isomorphism if and only if $\alpha$ itself is an isomorphism. If $\alpha \neq \beta$, then the given TF-isomophism $(\alpha,\beta)$ is essentially different from a usual isomorphism and hence we call $(\alpha,\beta)$ a \textit{non-trivial} TF-\textit{isomorphism}.  If $(\alpha,\beta)$ is a non-trivial TF-isomorphism from a mixed graph $G$ to a mixed graph $H$, the bijections $\alpha$ and $\beta$ need not necessarily be isomorphisms from $G$ to $H$. This is illustrated by the graphs in Figure \ref{fig:tfisomorph01}, 
examples found in \cite{lms2}, and also others presented below. \\

When $G=H$, $(\alpha,\beta)$ is said to be a TF-\textit{automorphism} and it is again called non-trivial if $\alpha \neq \beta$. The set of all TF-automorphisms of $G$ with multiplication defined by $(\alpha,\beta)(\gamma,\delta) = (\alpha \gamma, \beta \delta)$ is a subgroup of $S_{V(G)} \times S_{V(G) }$ and it is called the \textit{two-fold automorphism group} of $G$ and is denoted by \auttfg.  Note that if we identify an automorphism $\alpha$ with the TF-automorphism $(\alpha,\alpha)$, then Aut$(G) \subseteq$ \auttfg. When a graph has no non-trivial TF-automorphisms, Aut$(G)= $\auttfg.  It is possible for an asymmetric graph $G$, that is a graph with $|$Aut$(G)| = 1$, to have non-trivial TF-automorphisms. This was one of our main results in \cite{lms2}. \\

\bigskip

The main theme of this paper is \emph{stability of graphs}, an idea introduced by Maru{\v{s}}i{\v{c} et al. \cite{Scapsalvi1} and studied extensively by others, most notably by Wilson \cite{Wilson01} and Surowski \cite{Surowski1}, \cite{Surowski2}. Let $G$ be a graph and let $\BG$ be its \emph{canonical double cover} or \emph{duplex}. This means that V$(\BG)$ $=$ V$(G) \times \ZZ _{2}$ and  if $\{(u,0),(v,1)\}$ and $\{(u,1),(v,0)\}$ are edges of $\BG$ if and only if $\{u,v\}$ is an edge of $G$. One may think of the second entry in the notation used for vertices of $\BG$, that is $0$ or $1$ as colours. Recall that the graph $\BG$ is bipartite and we may denote its colour classes by $V_{0}= V \times \{0\}$ and $V_{1}= V \times \{1\}$ containing vertices of the type $(u,0)$ and $(u,1)$ respectively. A graph is said to be unstable if $\aut(G) \times \ZZ_{2}$ is a proper subgroup of $\aut \ \BG$. The elements of $\aut \ \BG$ $\setminus$ $\aut(G) \times \ZZ_{2}$ will be called \emph{unexpected automorphisms} of $\BG$. In other words, a graph $G$ is unstable if at least one element of $\aut \ \BG$ is not a lifting of some element of $\aut \ G$. In this paper, we shall investigate the relationship between the stability of the graph $G$ and its two-fold automorphism group \auttfg.


\section{Unstable Graphs and TF-automorphisms}

   
Consider $\aut \ \BG$. Let $\Sigma$ be the set-wise stabiliser of $V_{0}$  in $\aut \ \BG$, which of course coincides with the set-wise stabiliser of $V_{1}$. Note that every $\sigma \in \Sigma$ also fixes $V_{1}$ set-wise. We will show that it is the structure of $\Sigma$ which essentially determines whether $\BG$ has {\it{unexpected automorphisms}} which cannot be lifted from automorphisms of $G$.  The following result, which is based on the first result in \cite{Scapsalvi2}, Lemma 2.1 implies that these unexpected automorphisms of $\BG$ arise if the action of $\sigma$ on $V_{0}$ is not mirrored by its action of $V_{1}$.

{\Lem{Let $ f : \Sigma \rightarrow \mbox{\it{Sym}}(V) \times \mbox{\it{Sym}}(V)$ be defined by $ f : \sigma \mapsto (\alpha,\beta)$ where $\alpha$ is such that $(\alpha(v),0) = \sigma (v,0)$ and $(\beta(v),1) = \sigma(v,1)$, that is $\alpha$, $\beta$ extract from $\sigma$ its action on $V_{0}$ and $V_{1}$ respectively. Then: \begin{enumerate}
\item{$ f$ is a group homomorphism;}
\item{$ f$ is injective and therefore $ f: \Sigma \rightarrow  f(\Sigma)$ is a group automorphism; }
\item{$ f (\Sigma) = \{ (\alpha,\beta) \in Sym(V) \times Sym (V): x$ is adjacent to $y$ in $G$ if and only if $\alpha(x)$ is adjacent to $\beta(y)$ in $G\} $ that is, $ f(\Sigma)$= \auttfg, that is, $(\alpha,\beta)$ (the ordered pair of separate actions of $\sigma$ on the two classes) is a TF-automorphism of $G$. }
\end{enumerate}  
 }\label{lem:stabilitytfmain01}}

 {\proof{The fact that $ f$ is a group homomorphism, that is,  that $ f$ is a structure preserving map from $\Sigma$ to $\mbox{\it{Sym}}(V)\times \mbox{\it{Sym}}(V)$  follows immediately from the definition since for any $\sigma_{1}$, $\sigma_{2} \in \Sigma$ where $ f(\sigma_{1}) = (\alpha_{1},\beta_{1})$ and $ f(\sigma_{2} )= (\alpha_{2},\beta_{2})$,  $ f (\sigma_{1})  f (\sigma_{2})$ $= (\alpha_{1}\beta_{1})(\alpha_{2}\beta_{2}) = (\alpha_{1}\alpha_{2},\beta_{1}\beta_{2}) =  f (\sigma_{1}\sigma_{2})$.  This map is clearly injective and therefore $ f: \Sigma \rightarrow  f(\Sigma)$ is a group automorphism. \\
 
 \noindent Consider an arc $((u,0),(v,1))$ and then note that since $\sigma \in \Sigma \subseteq \aut\ \BG$, $(\sigma(u),0),(\sigma(v),1))$ is also an arc of $\BG$.  By definition, this arc may be denoted by $(\alpha(u),0), \beta(v),1)$ and, following the definition of $\BG$, it exists if and only if $(\alpha(u),\beta(v))$ is an arc of $G$. Hence $ f$ maps elements of $\Sigma$ to $(\alpha,\beta)$ which clearly take arcs of $G$ to arcs of $G$. This implies that $(\alpha,\beta)$ is a TF-automorphism of $G$ and hence $ f(\Sigma) =$\auttfg.
 }\lesta}
   
  \bigskip

  
  As shown in \cite{lms2}, Proposition 3.1, if $(\alpha,\beta)\in Aut^{TF} G$ then $(\gamma,\gamma^{-1})\in Aut^TF G$ where $\gamma=\alpha\beta^{-1})$. This means that for any edge $\{x,y\}$ of $G$, $\{\gamma(x),\gamma^{-1}(y)\}$ is also
an edge. A permutation $\gamma$ of $V(G)$ with this property is called an \emph{anti-automorphism}. Such maps possess intriguing applications to the study of cancellation of graphs in direct products with arbitrary bipartite graphs, that is, the characterisation of those graphs $G$ for which $G\times C \simeq H\times C$ implies $G\simeq H$, whenever $C$ is a bipartite graph (see \cite{HB} Chapter 9). The second part of Theorem 3.2 could be rephrased as follows: ``$G$ is unstable if and only if it has an anti-automorphism of order different from 2". Note that the existence of an anti-automorphisms of order 2 does not imply instability since such a map corresponds to a  trivial TF-automorphisms.


 {\Thm{Let $G$ be a graph. Then $\aut \ \BG = \mbox{\auttfg} \rtimes \ZZ_{2}$. Furthermore, $G$ is unstable if and only if it has a non-trivial \it{TF}-automorphism. }\label{thm:stabilitytfmain01b}}

 {\proof{ From Lemma \ref{lem:stabilitytfmain01}, $f(\Sigma) = $ \auttfg $\ $ which must have index 2 in $\aut \ \BG$. The permutation $\delta(v,\varepsilon) \mapsto (v,\varepsilon+1)$ is an automorphism of $\BG$  and $\delta \not \in f(\Sigma)$. Then $\aut \BG$ is generated by $f(\Sigma)$ and $\delta$. Furthermore, $f(\Sigma) \cap \langle\delta\rangle = \id$ and $f(\Sigma) \lhd \aut \ \BG$ being of index 2.\\
 
 Since $\aut\  \BG = \mbox{\auttfg} \rtimes \ZZ_{2}$, $G$ is stable if and only if \auttfg$\ = \aut\ G$. 
     }\lesta}

 \begin{figure}[h]
\centering
\includegraphics [width=9 cm, height= 4.5 cm] {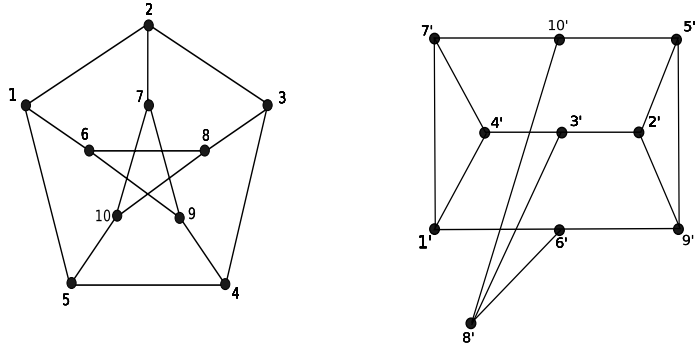}
\caption{A stable and unstable graph which are TF-isomorphic.}
\label{fig:altcon011}
\end{figure}

It is natural to ask whether it can happen  that a stable graph is TF-isomorphic to an  unstable one. The answer is yes and an example is shown in Figure \ref{fig:altcon011} with the Petersen graph being stable and the other graph which is TF-isomorphic to it being unstable. Both graphs have the same bipartite canonical double cover since they are TF-isomorphic. The reader may refer to \cite{porcu} to learn more about graphs having the same canonical double cover. We should point out here (as noted by Surowski in \cite{Surowski1}) that if a graph $G$ is stable, that is, $\aut \BG = $ $\mbox{\aut G} \rtimes \ZZ_{2}$, then the semi-direct product must be a direct product because $\aut G$ is normal in $\aut \BG$ since it has index $2$ and also $\ZZ_{2}$ is normal since its generator commutes with every element of $\aut \BG$, by stability. In fact, as Surowski comments, the stability of $G$ is equivalent to the centrality of $\ZZ_{2}$ in $\aut \BG$, which is the lift of the identity in $\aut G$.\\

\begin{figure}[h]
\centering
\includegraphics [width=9.5cm, height= 5.7 cm] {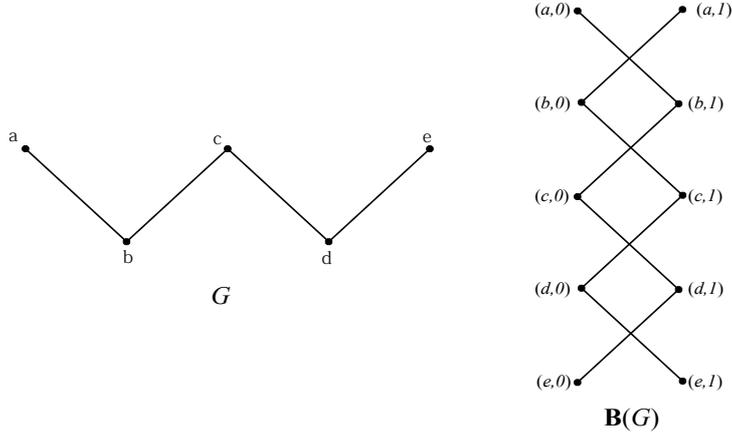}
\caption{An example used to show how a TF-automorphism  $(\alpha,\beta)$ of $G$ can be obtained from an automorphism $\sigma$ of $\BG$.}
\label{fig:instabilitytfmain01}
\end{figure}

 It is worth noting that the ideas explored in the proof of Lemma \ref{lem:stabilitytfmain01} may be used to extract TF-automorphisms of a graph $G$ from automorphisms of $\BG$ which fix the colour classes. In fact, let $\sigma$ be such an automorphism. Define the permutations $\alpha$ and $\beta$ of V$(G)$ as follows: $\alpha(x) = y$ if and only if $\sigma (x,0) = (y,0) $ and $\beta(x) = y$ if and only if $\sigma(x,1) = (y,1)$. Then $(\alpha,\beta)$ is a TF-automorphism of $G$. We remark that $\alpha$ and $\beta$ are not necessarily automorphisms of $G$ as we shall show in the example shown in Figure \ref{fig:instabilitytfmain01}. The automorphism $\sigma$ is chosen so that it fixes one component of $\BG$ whilst being an automorphism of the other component. In order to have a more concise representation, we denote vertices of $\BG$ of the form $(u,0)$, that is, elements of the colour class $V_{0}$ by $u_{0}$ and similarly denote  vertices of the form $(u,1)$ in $V_{1}$ by $u_{1}$. Using this notation, $\sigma = (a_{0})(b_{1})(c_{0})(d_{1})(e)(a_{1}\ e_{1})(b_{0}\ d_{0})(c_{1})$. The permutations $\alpha$ and $\beta$ of $G$ are extracted from $\sigma$ as described in the proof of Lemma \ref{lem:stabilitytfmain01}. For instance, to obtain $\alpha$, we restrict the action of $\sigma$ to the elements of $V_{0}$, that is, those vertices of the form $(v,0)$ or $v_{0}$ when using the new notation and then drop the subscript.  Similarly, the permutation $\beta$ is obtained from the action of $\sigma$ restricted to $V_{1}$. Therefore, $\alpha = (a)(c)(e)(b\ d)$ and $\beta = (b)(d)(a\ e)(c)$. Note that neither $\alpha$ nor $\beta$ is an automorphism of $G$, but $(\alpha,\beta)$ is a TF-automorphism of $G$ which in turn can be lifted to the unexpected automorphism $\sigma$ of $\BG$. This example illustrates Lemma \ref{lem:stabilitytfmain01} since the graph $G$ is unstable and has a non-trivial TF-automorphism. \\

 The result of Lemma \ref{lem:stabilitytfmain01} and the subsequent example lead us to other questions regarding the nature of the permutations $\alpha$ and $\beta$ which, as discussed in the preceding example given in Figure \ref{fig:instabilitytfmain01}, may not be automorphisms of $G$.  \\

If $(\alpha, \id)$ is a non-trivial TF-automorphism of a graph $G$, then $G$ is not vertex-determining. In fact, since $\alpha \neq \id$ then  $\alpha (u)  = v$ for some $u \neq v$ and the TF-automorphism $(\alpha,\id)$ fixes the neighbours of $u$ and takes $u$ to $v$. Hence $u$ and $v$ must have the same neighbourhood set, which implies that $G$ is not vertex-determining.\\

We shall use this idea to prove some results below. An alternative way of looking at this is to consider Lemma \ref{lem:stabilitytfmain01} and to note that a graph $G$ is stable if and only if given $\sigma(v,0) = (\alpha(v),0)$, there exists no $\beta \neq \alpha$ such that $(\beta(v),1))=\sigma(v,0)$. Hence, as implied by Theorem \ref{thm:stabilitytfmain01b} a graph $G$ is stable if and only if $f(\Sigma) \subseteq \Delta_{V}$ where $\Delta_{V}$ is the diagonal group of $(\alpha,\beta)$, $\alpha$, $\beta$ automorphisms of $G$, with $\alpha=\beta$.\\

 {\Prop{If $(\alpha,\beta)$ is a non-trivial TF-automorphism of a graph $G$ but $\alpha$ and $\beta$ are automorphisms of $G$, then $G$ is not vertex-determining.}\label{prop:alphabetaunworthy01}}
 
 {\proof{ Since $\alpha$ is an automorphism of $G$, then $(\alpha,\beta)$ is a TF-automorphism, so the group \auttfg $\ $ must also contain $(\alpha,\beta)(\alpha^{-1},\alpha^{-1})$ $=$ $(\id, \beta \alpha^{-1})$. Since $\alpha\beta^{-1} \neq \id$, let $u$ be a vertex such that $v = \alpha\beta^{-1}(u)$ is different from $u$. Then for each neighbour $w$ of $u$ the arc $(w,u)$ is taken to the arc $(w,v)$, so that $N(u)$ is contained in $N(v)$ and vice versa. Therefore $u$ and $v$ have the same neighbourhood, so $G$ is not vertex-determining.
 }\lesta}

 {\Prop{If $(\alpha,\beta)$ is a non-trivial TF-automorphism of a graph $G$ and $\alpha$, $\beta$ have a different order, then $G$ is not vertex-determining.}\label{prop:alphabetaunworthy02}}
 
 {\proof{Let $(\alpha,\beta)$ be an element of \auttfg $\ $with the orders of $\alpha$ and $\beta$ being $p$ and $q$ respectively and assume without loss of generality that $p < q$.
 Since \auttfg $\ $ is a group, $(\alpha,\beta)^{p} = (\alpha^{p},\beta^{p}) = (\id, \beta^{p})$ must also be in \auttfg. The same argument used in the proof of Proposition \ref{prop:alphabetaunworthy01} holds since $\beta^{p} \neq \id$. Hence $G$ is not vertex-determining.
 }\lesta}

 Proposition \ref{prop:alphabetaunworthy01} and Proposition \ref{prop:alphabetaunworthy02} are equivalent to their counterparts in \cite{Scapsalvi1} in which they are stated in terms of adjacency matrices. In \cite{Scapsalvi1}, it is shown that if a graph $G$ is unstable but vertex-determining and $(\alpha,\beta)$ is a non-trivial TF-automorphism of $G$, then $\alpha$ and $\beta$ must not be automorphisms of $G$ and must have the same order. This then gives us information about automorphisms $\sigma$ of $\BG$ which are liftings of TF-automorphisms of $G$.

  \section{Triangles}
 
In this section we shall study the behaviour of a non-trivial TF-automorphism of a graph $G$ acting on a subgraph of $G$ isomorphic to $K_{3}$ with the intent of obtaining information regarding the stability of graphs which have triangles as  a basic characteristic of their structure. Strongly regular graphs in which every pair of adjacent vertices have a common neighbour are an example. The stability of such graphs has been studied by Surowski \cite{Surowski1}.  We believe that this section is interesting because it is a source of simple examples of unstable graphs and also because a detailed analysis of what happens to triangles can throw more light on TF-automorphisms of graphs.\\

  \begin{figure}[h]
\centering
\includegraphics [width=9cm, height= 13.2 cm] {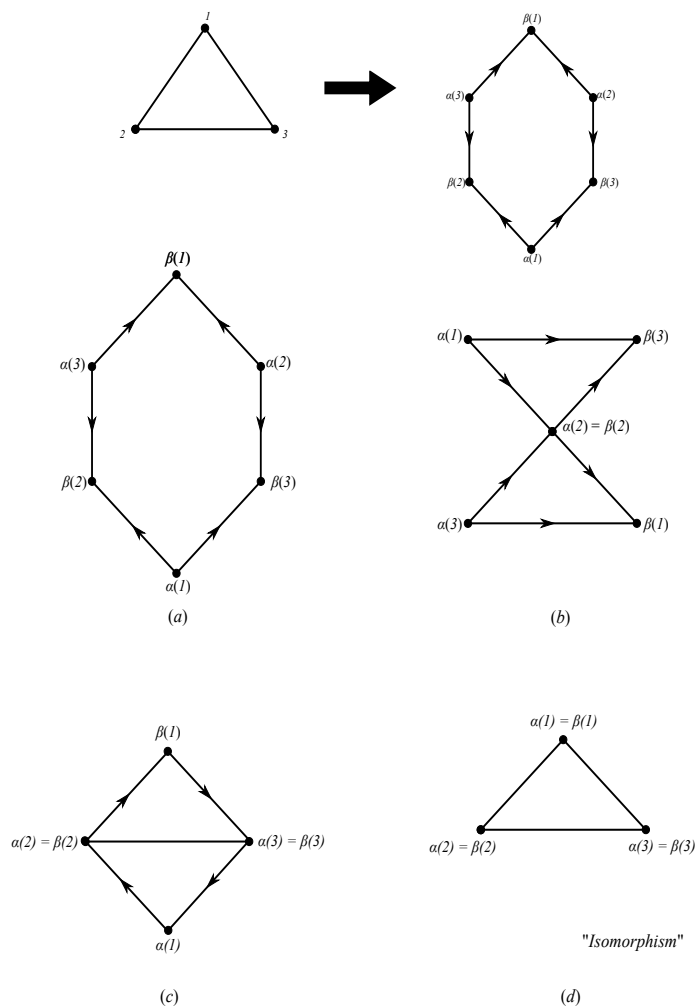}
\caption{The configurations of possible images of a triangle under the action of a non-trivial TF-automorphism as described in Proposition \ref{Prop:surowskitriangle01}.}
\label{fig:trianglec6and2tri}
\end{figure}

{\Prop{Let $(\alpha,\beta)$ be a TF-isomorphisms from a graph $G$ to a graph $G'$. The action of $(\alpha,\beta)$ on some subgraph $H \cong K_{3}$ yields either (a) a closed Z-trail of length 6 with no repeated vertices, (b) a pair of oriented Z-connected triangles with exactly one common vertex, (c) a pair of oriented triangles with exactly two common vertices or (d) an undirected triangle as illustrated in Figure \ref{fig:trianglec6and2tri}(a),(b),(c) and (d) respectively.}\label{Prop:surowskitriangle01}}

{\proof{Since $\alpha\neq \beta$, there exists some vertex $1$ such that $\alpha(1)\neq\beta(1)$.  Consider the triangle $H$ such that V$(H) = \{1,\ 2,\ 3\}$. The semi-closed Z-trails covering $H$ must be taken to the Z-trails $\alpha(1)\longrightarrow\beta(2)\longleftarrow\alpha(3)\longrightarrow \beta(1)$ and $\beta(1)\longleftarrow\alpha(2)\longrightarrow\beta(3)\longleftarrow \alpha(1)$ which together form a closed Z-trail of length 6. Vertices in Z-trails can be repeated and this has to be considered when studying the embedding of the given Z-trail within the graph. If no vertex is repeated then clearly we have a closed Z-trail of length 6 covering a subgraph isomorphic to $C_{6}$. If not, we have to consider all options. We are already assuming that $\alpha(1)\neq\beta(1)$.  Then, we can have $\alpha(2)=\beta(2)$ or $\alpha(3) = \beta(3)$ as in Figure \ref{fig:trianglec6and2tri}(b) or both as shown in Figure \ref{fig:trianglec6and2tri}(c) or, if $\alpha(1)=\beta(1)$, $\alpha(2)=\beta(2)$ and $\alpha(3)=\beta(3)$, we have an internal isomorphism as shown in Figure \ref{fig:trianglec6and2tri} (d).}\lesta}

In general, when a TF-automorphism $(\alpha,\beta)$ acts on the arcs of $G$ it maps any triangle $H$ into another triangle if $\alpha(x) = \beta(x)$ for every vertex $x$ of the triangle or it fits one of the other configurations described by Proposition \ref{Prop:surowskitriangle01} which are illustrated in Figure \ref{fig:trianglec6and2tri}.  If a graph $G$ in which every edge lies in a triangle is unstable, then it must have a non-trivial TF-automorphism which follows one of these configurations.\\


 \begin{figure}[h]
\centering
\includegraphics [width=7cm, height=  8cm] {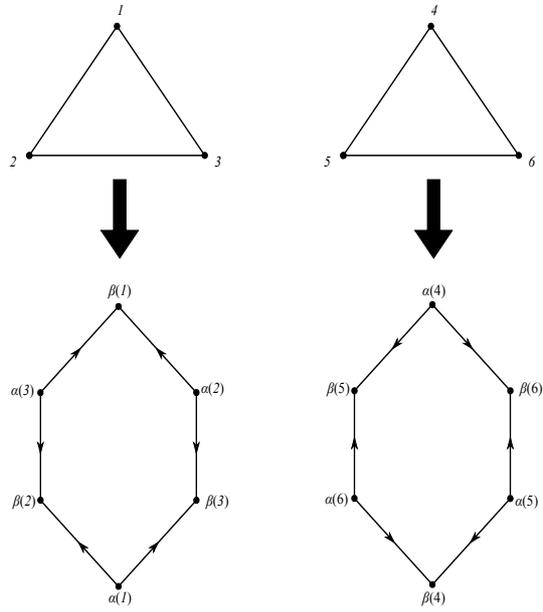}
\caption{Two triangles mapped to a $C_{6}$.}
\label{fig:2tri2tritfa}
\end{figure} 

\begin{figure}[h]
\centering
\includegraphics [width=7cm, height= 8.5 cm] {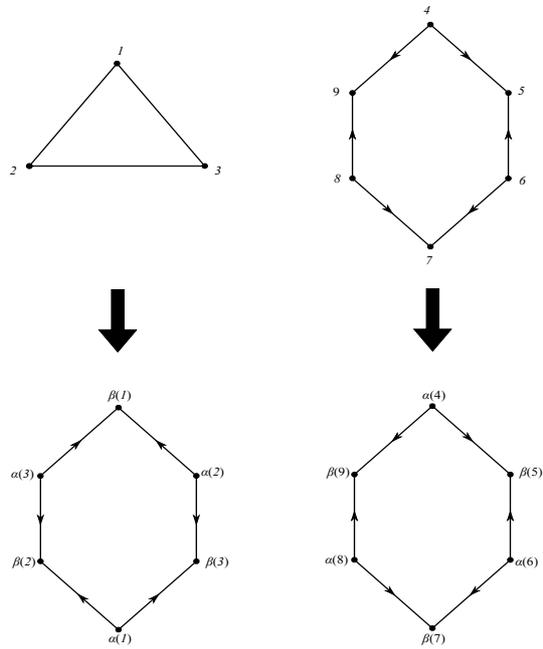}
\caption{A triangle and one closed Z-trail of length 6 covering the edges of a $C_{6}$ mapped to a $C_{6}$.}
\label{fig:2tric6tfa}
\end{figure} 

\clearpage

{\Prop{Let $(\alpha,\beta)$ be a TF-isomorphism from $G$ to $G'$. When the TF-isomorphism acting on $K_{3}$, a subgraph of $G$, yields a closed Z-trail of length 6 with no repeated vertices as shown in figure \ref{fig:trianglec6and2tri}(a), either two triangles with no common vertex or a triangle and a closed $Z$-trail of length $6$ are mapped to a subgraph isomorphic to $C_{6}$. In the cases when the TF-isomorphism acting on a $K_{3}$ yields the images illustrated in Figure \ref{fig:trianglec6and2tri} (b), (c), the pair of triangles which are either mapped to two triangles with exactly one common vertex or to two triangles with exactly one common edge must be connected. }\label{Prop:tfatfbgraph01}}

{\proof{Refer to Figure \ref{fig:trianglec6and2tri}. In the case illustrated in Figure \ref{fig:trianglec6and2tri}(a) the arcs of one closed Z-trail $P$ of length 6 can be the co-domain the arcs of a  triangle $H$. The Z-trail $P'$ obtained by reversing the arcs of $P$ can be the co-domain of another triangle $K$. We claim that $H$ and $K$ are vertex disjoint.  In fact, suppose not and assume that the two triangles have a common vertex $u$. The pair of vertices $\alpha(u)$ and $\beta(u)$ where $\alpha(u) \neq \beta(u)$ are in both $P$ and $P'$ and this is contradiction as the in-degree of $\alpha(u)$ and similarly the out-degree of $\beta(u)$ must be zero and this makes it impossible to identify arcs of $P$ with arcs of $P'$ to form the edges of a $C_{6}$. Figure \ref{fig:2tri2tritfa} shows an example where setting $\alpha(1) =\beta(5)$,  $\beta(2)=\alpha(6)$,  $\alpha(1)=\beta(4)$,  $\beta(3)=\alpha(5)$,  $\alpha(2)=\beta(6)$, $\beta(1) = \alpha(4)$ would be one way of associating one directed $C_{6}$ with the other so that the alternating connected circuits form an undirected $C_{6}$.
The other possibility is illustrated in Figure \ref{fig:2tric6tfa}. In this example $\beta(1)=\alpha(4)$, $\alpha(2) = \beta(5)$, $\beta(3)=\alpha(6)$, $\alpha(1)=\beta(7)$, $\beta(2)=\alpha(8)$ and $\alpha(3)=\beta(9)$ so that a closed Z-trail of length $6$ covering a $K_{3}$ is mapped to a closed Z-trail of length $6$ covering half of the arcs of an undirected $C_{6}$ whilst the rest of the arcs come from a Z-trail of length $6$ covering half of the arcs of another subgraph isomorphic to $C_{6}$. The $K_{3}$ and the $C_{6}$ in the domain of the TF-isomorphism cannot have a common vertex and the proof is analogous to the one concerning the former case. \\

The proof for the remaining cases may be carried over along the same lines as the above. Refer to Figure \ref{fig:trianglec6and2tri}. We observe that we that the Z-trail $P_{1}$ described by $\beta(3)\longleftarrow \alpha(2) \longrightarrow \beta(1)$, the Z-trail $P_{2}$ described by $\alpha(1) \longrightarrow \beta(2) \longleftarrow \alpha(3)$ and the Z-trails $P_{1}'$ and $P_{2}'$ obtained by the arcs of $P_{1}$ and $P_{2}$ respectively would imply by the conservation of Z-trails, that in the pre-image of the subgraph, there are four Z-trails passing through the vertex labelled 2. This can only be possible if the triangles in the pre-image have a common vertex.}\lesta}
\bigskip

 \begin{figure}[h]
\centering
\includegraphics [width=7cm, height= 7 cm] {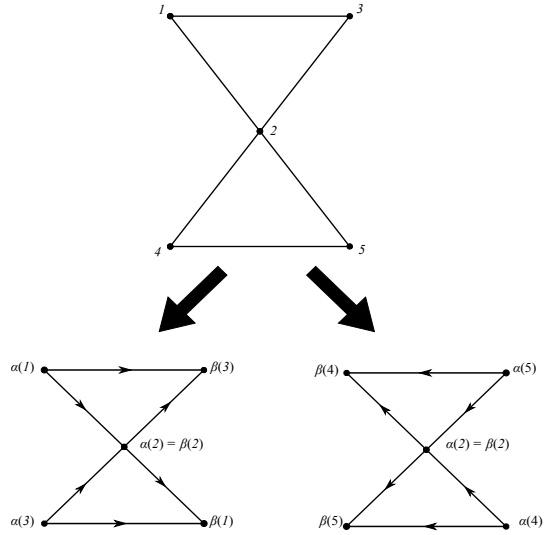}
\caption{An example to illustrate the result of Proposition \ref{Prop:tfatfbgraph01} where $\alpha(1)=\beta(4)$, $\beta(3)=\alpha(5)$, $\alpha(3) =\beta(5)$ and $\beta(1) = \alpha(4)$.}
\label{fig:2tri2tri}
\end{figure}
 
  \begin{figure}[h]
\centering
\includegraphics [width=7cm, height= 4 cm] {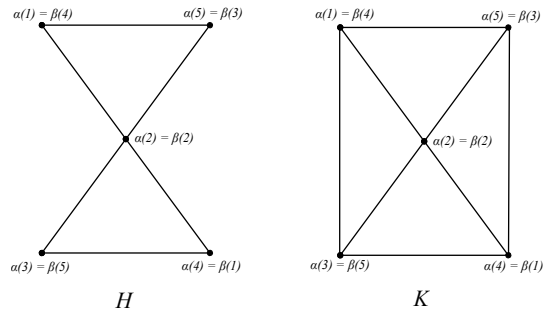}
\caption{The smallest unstable graphs where a triangle is taken to a two directed triangles sharing a vertex.}
\label{fig:2tri2tri02}
\end{figure}

\begin{figure}[h]
\centering
\includegraphics [width=8cm, height= 4 cm] {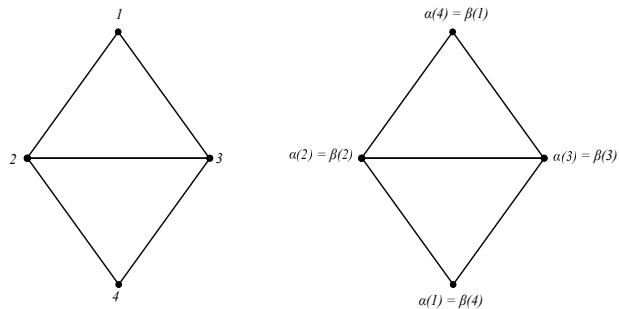}
\caption{The smallest unstable graph which has a TF-automorphism taking a triangle to the mixed graph illustrated in Figure \ref{fig:trianglec6and2tri}(c).}
\label{fig:2tri2tri02tf3}
\end{figure}

\clearpage

Figure \ref{fig:2tri2tri02} shows the smallest unstable graphs which have a TF-automorphism which takes a triangle to a pair of directed triangles with a common vertex as illustrated in figure \ref{fig:trianglec6and2tri}(b).  Figure \ref{fig:2tri2tri02tf3} shows the smallest graph which has a nontrivial TF-automorphism which maps a triangle to the mixed graph illustrated in Figure \ref{fig:trianglec6and2tri}(c).\\


\section{Unstable graphs of arbitrarily large diameter}

In this section, we present a method of constructing unstable graphs of an arbitrarily high diameter. \\

If $H$, $K$ are graphs, let $[H,K]$ be the graph whose vertex set is the union V$(H)\cup$ V$(K)$ and whose edge set
is the union of E$(H)$, E$(K)$ plus the edges of the complete bipartite graph with classes V$(H)$ and V$(K)$. More generally, if $H_0,\ H_1,\ ...,\ H_{m-1}$  are graphs (where the labels are in $\ZZ_{m}$, where $m>1$ is an integer), let $G=[H_0,\ H_1,\ ...,\ H_{m-1}]$ be the graph whose vertex set is the union of all V$(H_i)$ and whose edge set is the union of  all E$([H_i,H_{(i+1)}])$. In other words, $G$ contains all vertices and edges of the graphs $H_i$, plus all edges of the complete bipartite graph connecting two consecutive $H_i$'s.\\

Now, assume that none of the $H_{i}$ has isolated vertices. Let $(\alpha_i,\beta_i):H_i \rightarrow H_{i+1}$ be TF-isomorphisms as $i$ runs over $\mathbb{Z}_m$. Assume that the product 
\[(\alpha_0,\beta_0)(\alpha_1,\beta_1)...(\alpha_{m-1},\beta_{m-1}) = (\id,\id).\]
Note that the latter assumption is not a restriction, because one can always take $(\alpha_0,\beta_0)$ as the inverse of the product of the remaining TF-isomorphisms.\\

{\Thm{With the above assumptions, let $G$ $=$ $[H_0,\ H_1,\ ...,\ H_{m-1}]$. Define  two permutations $\alpha$, $\beta$ of V$(G)$ as
follows. For $v\in$ V$(G)$, let $i$ be such that $v\in$ V$(H_i)$; then set $\alpha(v)=\alpha_{i}(v)$ and $\beta(v)=\beta_{i}(v)$. \\

Then the following hold:\begin{enumerate}

\item{$(\alpha,\beta)$ is a TF-automorphisms of $G$;\\}
\item{diam $G$ $=k$ $= (m+e)/2$, where $e=0$ if $m$ is even and $e=1$ if $m$ is odd;\\}
\item{Each edge of $G$ belongs to a triangle;\\}
\item{Letting $m=2k+e$, whenever $d(v,w)=k-1$, there is a vertex $u$ adjacent to $w$ such that d$(v,u)=k$.}
\end{enumerate}}\label{Thm:surowskicounter01}}

{\proof{ First note that if $uv$ is an arc of $G$ and both $u$, $v$ belong to the same $H_i$, then the image of $uv$ is an arc of $H_{i+1}$, hence of $G$,
because $(\alpha,\beta)$ acts like $(\alpha_i,\beta_i)$ in $H_i$. If $u$, $v$ do not belong to the same $H_i$, then they belong to consecutive graphs,
say $H_i$, $H_{i+1}$, so $\alpha(u)$ and $\beta(u)$ belong to the consecutive graphs $H_{i+1}$, $H_{i+2}$, and are adjacent because all
the arcs between these two graphs belong to $G$. This proves (1).\\

\noindent Concerning distance, a path from $u$ in $H_0$ to $v$ in $H_s$, where $s=(m-e)/2$, must pass through all graphs $H_1,\  H_2,\ ...,\ H_{s-1}$
or else $H_{m-1},\ H_{m-2},\ ...,\ H_{s+1}$. Since such a path can be found, d$(u,v)=s$. For two vertices in, say, $H_i$ and $H_j$
with $i\neq j$, the same argument shows that d$(u,v)$ cannot exceed $s$. Finally, if $u,\ v$ lie in the same $H_i$, they have a
common neighbour in $H_{i+1}$, then d$(u,v)$ is less or equal than $2$ (regardless to their distance within $H_i$). This proves (2).\\

\noindent What about triangles?  If $uv$ is an edge of some $H_i$ then letting $w\in H_{i+1}$ the vertex $w$ is adjacent to both $u$ and $v$, hence
$uvw$ is a triangle. If $uv$ is an edge of some $[H_i,\  H_{i+1}]$, say with $u \in H_i$ and $v \in H_{i+1}$, take any neighbour $w$ of $u$ in $H_i$
(recall the assumption about no isolated vertices) and get the triangle $uvw$. This proves (3).\\

\noindent Finally, since $m=2k+e$, diam $G$ is at least $k+e$. Let us remark that if d$(v,w)=k-1$ and, 
say, $v\in$ V$(H_i)$, then $w\in$ V$(H_j)$ where $j$ is either $i+2$ or $i-2$. Take $u$ in $H_{i+3}$ or $H_{i-3}$, respectively. Then d$(v,u)=k$. 
This proves (4). Note that the structure of the single $H_i$'s is immaterial here, because two vertices at distance at least $k-1$ from each other must come from a different $H_{i}$.}\lesta}

Until now, we did not mention that the concerned TF-isomorphisms are non-trivial, so all the above would work fine for
the case of isomorphisms too. But adding the hypothesis that at least one of them is non-trivial, the obtained graph $G$
has a non-trivial TF-isomorphism, namely $(\alpha,\beta)$ as described above. The statements (1)-(3) show that there are
unstable graphs of arbitrarily high diameter, where each edge belongs to a triangle. \\

 \begin{figure}[h]
\centering
\includegraphics [width=8.4 cm, height= 12cm] {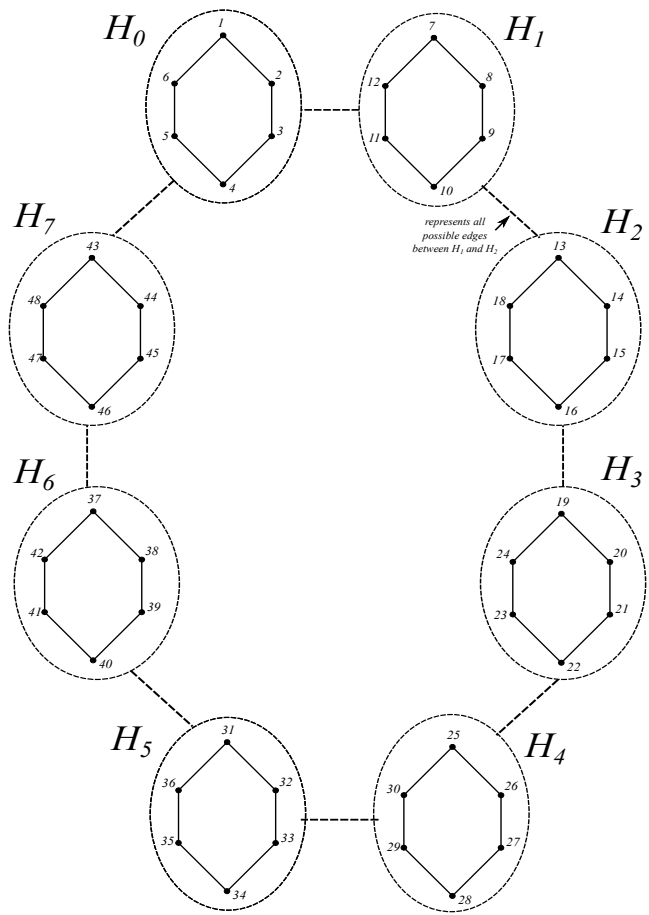}
\caption{A graph constructed using Theorem \ref{Thm:surowskicounter01}.}
\label{fig:surowski_counter01a}
\end{figure}

Surowski \cite{Surowski1}, \cite{Surowski2} proved various results concerning graph stability. In \cite{Surowski1} Proposition 2.1, he claims that if $G$ is a connected graph of diameter $d \geq 4$ in which every edge lies in a triangle, then $G$ is stable. However, by taking $m \geq7$ in Theorem \ref{Thm:surowskicounter01} we get infinitely many counterexamples to this claim by taking all the $H_{i}$ isomorphic to the same vertex-determining bipartite graph, because such a vertex-transitive graph is unstable, therefore one can find a non-trivial TF-isomorphism from $H_{i}$ to $H_{i+1}$, and since these $H_{i}$ are isomorphic, the resulting graph $G$ is vertex-determining, and has diameter $k\geq 4$  and unstable.

 One of these counterexamples is illustrated in Figure \ref{fig:surowski_counter01a}.   \\

We detected one possible flaw in Surowski's proof. It is claimed in \cite{Surowski1} that whenever an automorphism of $\BG$ fixes $(v, 1)$ it also fixes $(v,-1)$. We have not seen a proof of this result. Besides, in our last example, $G$ has a non-trivial fixed-point-free TF-automorphism, which implies that $\BG$ has a fixed-point-free automorphism that that fixes the colour classes. This claim is also used in \cite{Surowski1} Proposition 2.2 which states that if $G$ is a strongly regular graph with $k > \mu \not = \lambda \geq 1$, then $G$ is stable. Hence, we believe that at this point, the stability of strongly regular graphs with these parameters requires further investigation. \\

\section{Concluding Remarks}

The use of TF-isomorphisms in the study of stability of graphs provides a fresh outlook which allows us to view facts within a more concrete framework
and also provides tools to obtain new results. For instance, we can investigate the structure of the given graph without actually requiring
to lift the graph to its canonical double cover, but only having to reason within the original graph. Furthermore, the insights that we already have about
TF-isomorphisms of graphs may be considered to be new tools added to a limited toolkit. In particular, let us mention the idea of graph invariants under
the action of TF-isomorphisms, such as Z-trails, a topic which we have started to study in \cite{lms3}. To be able to find out how the subgraphs of a graph are related to other subgraphs within the graph itself in the case of unstable graphs fills a gap in our understanding of graph stability and using TF-isomorphisms appears to be a promising approach in this sense. We believe that this paper substantiates these claims. Furthermore, it motivates us to carry out further investigations. Some pending questions such as those concerning the stability of certain strongly regular graphs have already been indicated. The study of how TF-isomorphisms act on common subgraphs such as triangles is another useful lead. Nevertheless, the more ambitious aim would be the classification of unstable graphs in terms of the types of TF-automorphisms which they admit.

\nocite{lauri2} 
 
\nocite{Scapsalvi2}

\bibliography{reference.bib}
\bibliographystyle{plain}

\end{document}